# TRANSCENDENTAL FUNCTIONS WITH A COMPLEX TWIST

**Michael Warren**
**Department of Mathematics**
**Tarleton State University**
**Box T-0470**
**Stephenville, TX 76402**
**mwarren@tarleton.edu**

**Dr. John Gresham**
**Department of Mathematics**
**Tarleton State University**
**Box T-0470**
**Stephenville, TX 76402**
**jgresham@tarleton.edu**

**Dr. Bryant Wyatt**
**Department of Mathematics**
**Tarleton State University**
**Box T-0470**
**Stephenville, TX 76402**
**wyatt@tarleton.edu**

Abstract

In our previous paper, Real Polynomials with a Complex Twist [see http://archives.math.utk.edu/ICTCM/VOL28/A040/paper.pdf], we used advancements in computer graphics that allow us to easily illustrate more complete graphs of polynomial functions that are still accessible to students of many different levels. In this paper we examine the 3D graphical representations of selected transcendental functions over subsets of the complex plane for which the functions are real-valued. We visualize and find connections between circular trigonometric functions and hyperbolic functions.

## Introduction

In Warren, Gresham & Wyatt (2016), the authors established that students could benefit from a three-dimensional visualization of polynomial functions with real coefficients whose domain contains all complex inputs that produce real-valued outputs. Connections between the symbolic and graphical representations can stimulate student comprehension of functions (Lipp, 1994; Presmeg, 2014; Rittle-Johnson & Star, 2009). Students need to understand the concepts of domain, range, and the graphical representation of functions in order to be prepared for advanced mathematical study (Martinez-Planell & Gaisman 2012). The process of domain restriction involved in creating three-dimensional graphical representations of commonly studied functions provides a context for students to explore all these ideas.

When working with polynomial functions, instructors often introduce the existence of complex roots via the Fundamental Theorem of Algebra followed by symbolic methods for finding those roots. Traditional graphical representations do not always support the solutions obtained symbolically. Domain restriction coupled with interactive 3D graphics in programs like GeoGebra help students to make the essential connections between function representations. For a detailed exploration of this process for polynomial functions along with instructions on creating the 3D polynomial curves in GeoGebra, see Warren, Gresham & Wyatt (2016) at http://archives.math.utk.edu/ICTCM/VOL28/A040/paper.pdf.

Unlike polynomial functions, it is rare for students to study the non-real roots of transcendental functions outside of a complex variables course. In preparation for advanced mathematical study, students could investigate some of the non-real behavior of transcendental functions through a three-dimensional graphical representation. Within the context of roots for transcendental equations, students can connect the concept of no real roots with the existence of non-real roots. Utilizing the process of domain restriction and the 3D graphing capabilities of GeoGebra, instructors can give students the tools to understand some of the complex behavior of transcendental functions.

## Trigonometric and Hyperbolic Functions

Consider the function $\sin x + 2$. The authors ask two simple questions, which inputs generate an output of two and which inputs generate an output of zero? The first question and its solution set is common in a Precalculus or Trigonometry course. The number of solutions is infinite and the solution set has symbolic representation $\{x | x = n\pi\}$ where $n \in \mathbb{Z}$. Students can verify the answers to this question by examining the graph, see *figure 1*. The symbolic and traditional graphical representations are well-aligned and both serve to strengthen student understanding of the function. The second question is also common but the typical solution set needs explanation. An equation such as $\sin x + 2 = 0$ is said to have no solutions because the domain is restricted to be all real numbers. This leads to a false understanding of the sine function that influences student understanding throughout further mathematics study. The traditional graphical representation supports this misunderstanding, see *figure 1*.

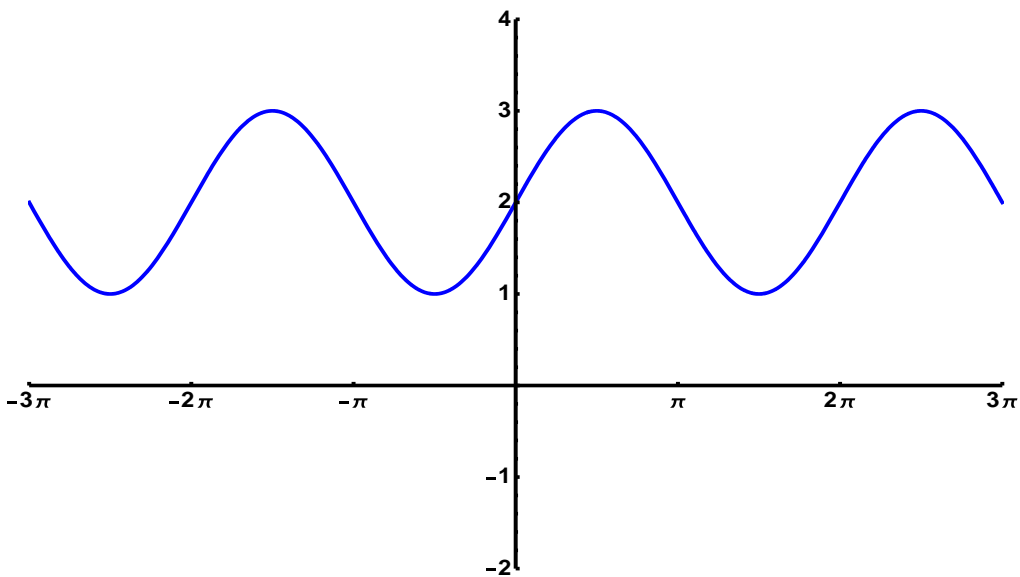

*Figure 1.* Graph of $\sin x + 2$. The solutions for $\sin x + 2 = 2$ are easily verified. The existence of solutions for $\sin x + 2 = 0$ are not apparent.

Domain restriction can be used to produce a three-dimensional visualization of the sine function. This visualization shows all input-output pairs where each input is a complex numbers and each output is a real number. All input-output pairs form a graphical representation of the function that can be visualized in a three-dimensional coordinate system. The three-dimensional visualization of the sine function in *figure 2* can be generated as follows.

Let $f(z) = \sin z$
where $z = x + iy$ for $x, y \in \mathbb{R}$ and $i^2 = -1$.
Then, $f(z) = \sin(x + iy) =$

(1) $\sin x \cosh y + i \cos x \sinh y$.

In order to produce a real output, $\cos x \sinh y = 0$. So either $\cos x = 0$, which occurs when $x = (2n+1)\frac{\pi}{2}$ where $n \in \mathbb{Z}$, or $\sinh y = 0$, which occurs when $y = 0$. The restricted domain is therefore given by,

(2) $y = 0$
and
(3) $x = (2n+1)\frac{\pi}{2}$.

When $y = 0$, the typical Cartesian graph of $\sin x$ is generated. When $x$ is an odd multiple of $\frac{\pi}{2}$, it appears that a hyperbolic cosine curve exists in that plane intersecting the real sine curve at its maximum or minimum. These hyperbolic cosine curves are obtained by substituting $x = (2n+1)\frac{\pi}{2}$ into (1):

$$f(x+iy) = \sin x \cosh y + i \cos x \sinh y = \begin{cases} \cosh y & \text{if } n \text{ is even}, x = (2n+1)\frac{\pi}{2} \\ -\cosh y & \text{if } n \text{ is odd}, x = (2n+1)\frac{\pi}{2} \end{cases}$$

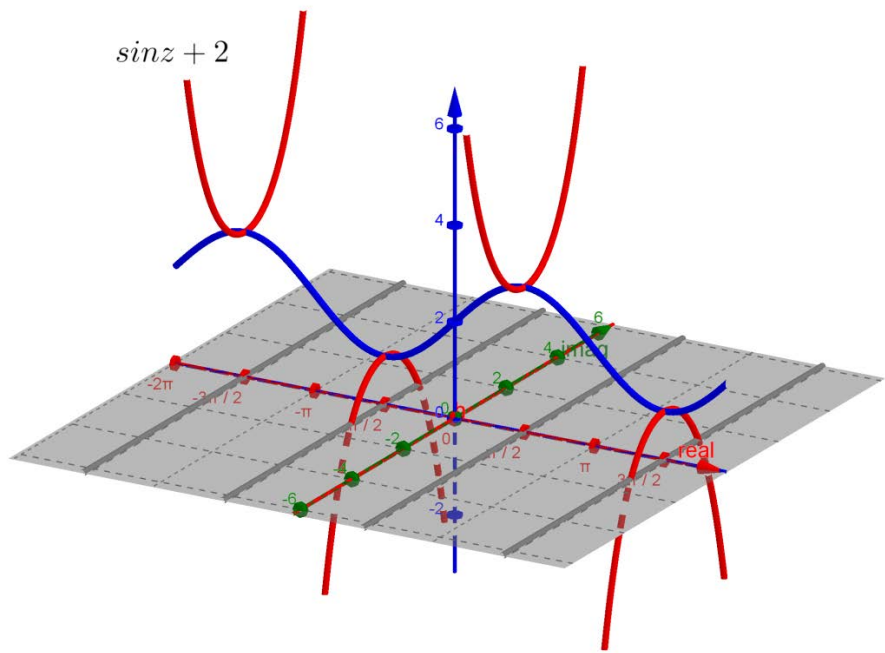

*Figure 2.* Graph of $\sin z + 2$. The existence of solutions for $\sin x + 2 = 0$ are and their values are easily verified.

This three-dimensional graphical representation allows a more complete understanding of the sine function. There are answers to the question about the zeros of the function $\sin z + 2$ but they are clearly all non-real and exist in conjugate pairs. The graphical representation gives meaning to the statement no real solutions that students and instructors often take for granted in trigonometric equations. It is interesting to note that the non-real curves for $\sin z$ are not trigonometric functions but hyperbolic cosines.

Next consider the secant function. Let $f(z) = \sec z = \sec(x + iy) =$

(4) $\dfrac{1}{\cos x \cosh y + i \sin x \sinh y}$.

In order to produce a real output, $\sin x \sinh y = 0$ and $\cos x \cosh y \neq 0$. Either $\sin x = 0$, which occurs when $x = n\pi$ where $n \in \mathbb{Z}$, or $\sinh y = 0$, which occurs when $y = 0$. Note that none of these values make $\cos x \cosh y = 0$. The restricted domain is therefore given by

(5) $y = 0$
and
(6) $x = n\pi$.

When $y = 0$, once again the typical Cartesian graph of $\dfrac{1}{\cos x} = \sec x$ is generated. When $x$ is an integer multiple of $\pi$, it appears that a hyperbolic secant curve exists in that plane intersecting the real secant curve at its maximum or minimum. The symbolic form of the non-real hyperbolic secant curves are obtained by substituting $x = n\pi$ into (4):

$$f(x + iy) = \frac{1}{\cos x \cosh y + i \sin x \sinh y} = \begin{cases} \operatorname{sech} y & if\ n\ is\ even, x = n\pi \\ -\operatorname{sech} y & if\ n\ is\ odd, x = n\pi \end{cases}$$

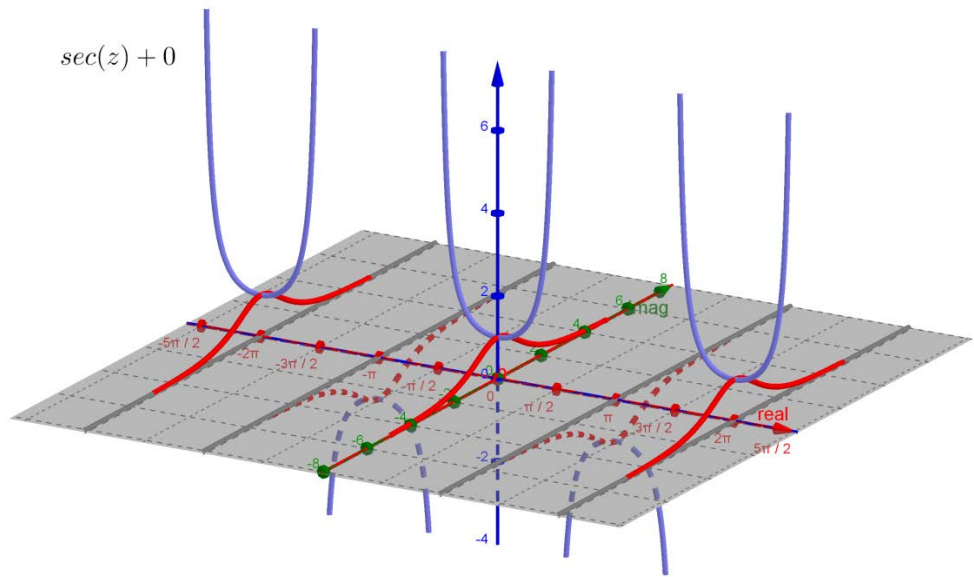

*Figure 3.* 3D Graph of $\sec z$. The non-real curves in red manifest as $\pm \operatorname{sech} y$ periodically when $x = n\pi$.

*Figure 3* illustrates the three-dimensional secant function. Note that the non-real hyperbolic secant curves exist between the local maximum and minimum values of the real secant curve. This behavior is similar to that observed in the sine function where the non-real hyperbolic cosine curves exist above and below the local maximum and minimum values respectively. In a sense, the non-real curves "fill in" the outputs missing from the real curves. Also similar to the sine function, the non-real curves for $\sec z$ are hyperbolic functions.

Since the sine and secant functions produce non-real hyperbolic curves, it is natural to follow with an inspection of a typically real-valued hyperbolic function. Let $f(z) = \cosh z = \cosh(x + iy) =$

(7) $\cosh x \cos y + i \sinh x \sin y.$

In order to produce a real output, $\sinh x \sin y = 0$. Either $\sinh x = 0$, which occurs when $x = 0$, or $\sin y = 0$, which occurs when $y = n\pi$ where $n \in \mathbb{Z}$. The restricted domain is therefore given by

(8) $x = 0$
and
(9) $y = n\pi$.

When $y = 0$, the typical Cartesian graph of $\cosh x$ is generated. When $y$ is any other integer multiple of $\pi$, another hyperbolic cosine curve is generated in that plane. When $x = 0$, there exists a cosine curve in the $yz$-plane where each hyperbolic cosine curve intersects at the local maximum and minimum values. *Figure 4* illustrates these curves. Note that all the output values not included in the range of the real hyperbolic cosine are included in the range of the three-dimensional hyperbolic cosine function. This behavior is consistent with the sine and secant functions examined earlier, with the exception of a zero output for the secant.

It is left to the reader to investigate the remaining trigonometric and hyperbolic functions.

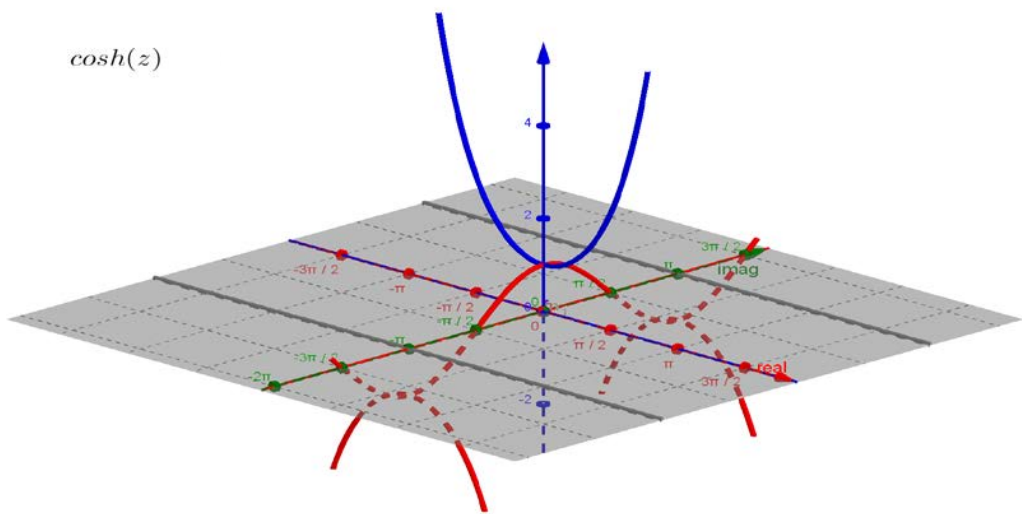

*Figure 4.* 3D graph of $\cosh z$. The graph of $\cos y$ is produced in the $yz$-plane with hyperbolic cosine curves intersecting at the local maximum and minimum values.

## The Exponential Function

Another well studied transcendental function worthy of examination is the exponential function, $e^z$. Consider the well-known equation $e^z + 1 = 0$. The traditional graphical representation, when the domain is restricted to all real numbers, leads students to believe that such an equation has no solutions. To combat this misconception, students can produce and explore the three-dimensional graphical representation of the exponential function.

Let $f(z) = e^z = e^{x+iy} =$

(10) $\qquad e^x(\cos y + i \sin y)$

In order to produce a real output, $e^x \sin y = 0$ which occurs when $y = n\pi$ where $n \in \mathbb{Z}$. When $y = 0$, the typical Cartesian graph of $e^x$ is generated. When $y$ is any other integer multiple of $\pi$, another exponential curve is generated in that plane. By substituting $y = n\pi$ where $n \in \mathbb{Z}$, (10) becomes

$f(z) = e^x \cos y$, where $y \in \{n\pi | n \in \mathbb{Z}\}$.

This produces the real and non-real curves shown in *figure 5*. Note that the exponential curves occur with a period of $2\pi$. The three-dimensional representation of the exponential function allows students to visualize the periodic behavior of the exponential function and gives meaning to the statement that there are no real solutions.

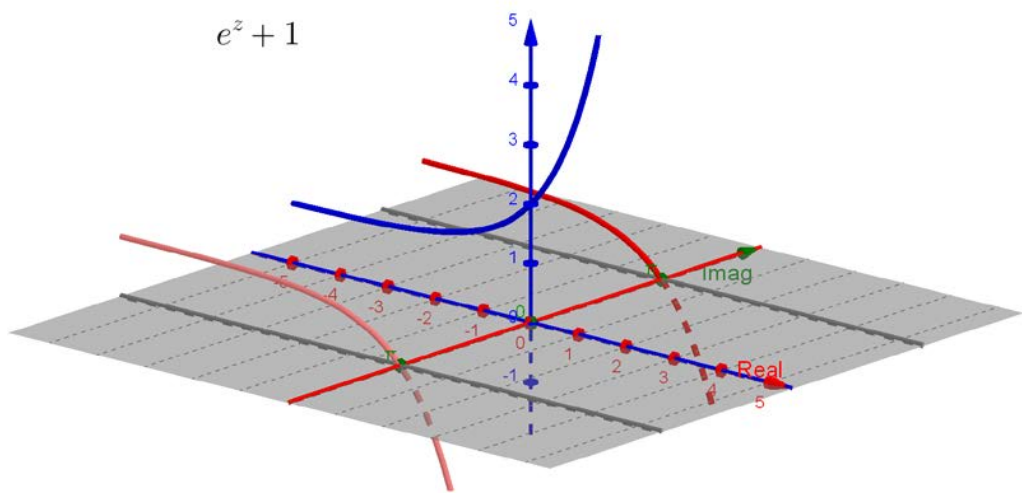

*Figure 5.* 3D graph of $e^z + 1$

## Conclusion

Students need to see the connections between the symbolic and graphical representations of a function. Traditional graphs of transcendental functions do not reveal the assumptions associated with claims that the equations $\sin x + 2 = 0$ and $e^x + 1 = 0$ have no solutions. Students fail to grasp the importance and necessity of stating there are no real solutions when teachers do not allow them to explore the non-real behavior of functions. With the technology available today, this ignorance is not necessary. Domain restriction, coupled with interactive 3D graphics in programs like GeoGebra, provide an excellent laboratory for discovery. Students now have the tools to investigate a more complete transcendental function with symbolic and graphical representations that support each other.

# GeoGebra Links

The sine function illustrator used in this paper can be found at https://ggbm.at/u8fHcJdk.

The secant function illustrator used in this paper can be found at https://ggbm.at/Mgc7CEMc.

The hyperbolic cosine function illustrator used in this paper can be found at https://ggbm.at/S8ZXU8bd.

The exponential function illustrator used in this paper can be found at https://ggbm.at/HD2JW5nZ.